 \documentclass[11pt]{article}
 \usepackage[top=1in,bottom=1in,left=1.5in,right=1.5in]{geometry}
  \usepackage{amsmath,amssymb}
  \usepackage{latexsym}
 \usepackage[dvips]{pstricks} 
  \usepackage{pst-node}
  \usepackage[dvips]{graphicx}

  \newcommand{\C}{\mathbb{C}}
  \newcommand{\F}{\mathbb{F}}

  \newcommand{\R}{\mathbb{R}}
  \newcommand{\Z}{\mathbf{Z}}

  \newcommand{\e}{\mathbf{e}}
  \newcommand{\f}{\mathbf{f}}
  \newcommand{\g}{\mathbf{g}}

  \newcommand{\U}{\mathbf{U}}
  \renewcommand{\u}{\mathbf{u}}
  \renewcommand{\v}{\mathbf{v}}
  \newcommand{\V}{\mathbf{V}}
  \newcommand{\w}{\mathbf{w}}
  \newcommand{\W}{\mathbf{W}}
  \newcommand{\x}{\mathbf{x}}
  \newcommand{\X}{\mathbf{X}}
  \newcommand{\y}{\mathbf{y}}
  \newcommand{\Y}{\mathbf{Y}}
  \newcommand{\z}{\mathbf{z}}
  \newcommand{\0}{\mathbf{0}}

  \newcommand{\cA}{\mathcal{A}}
  \newcommand{\cB}{\mathcal{B}}
  \newcommand{\cC}{\mathcal{C}}

  \newcommand{\cR}{\mathcal{R}}
  \newcommand{\cS}{\mathcal{S}}
  \newcommand{\cT}{\mathcal{T}}
  \newcommand{\cU}{\mathcal{U}}
  \newcommand{\cV}{\mathcal{V}}
  
  \newcommand{\cX}{\mathcal{X}}
  \newcommand{\cY}{\mathcal{Y}}

  \newcommand{\lan}{\langle}
  \newcommand{\ran}{\rangle}
  \newcommand{\an}[1]{\lan#1\ran}
  
  \newcommand{\hs}{\hspace*{\parindent}}
  \newcommand{\proof}{\hs \textbf{Proof.\ }}

  \newcommand{\Gr}{\mathop{\mathrm{Gr}}\nolimits}
  \newcommand{\Fr}{\mathop{\mathrm{Fr}}\nolimits}
  \newcommand{\trans}{^\top}
  \newcommand{\qed}{\hspace*{\fill} $\Box$\\}
  \newcommand{\opt}{\mathop{\mathrm{opt}}\nolimits}

  \newcommand{\dist}{\mathrm{dist}}

  \newcommand{\rS}{\mathrm{S}}

  \newcommand{\rank}{\mathrm{rank}}

  \newtheorem{theo}{\bfseries \hs Theorem}[section]
  \newtheorem{defn}[theo]{\bfseries \hs Definition}
  \newtheorem{prop}[theo]{\bfseries \hs Proposition}
  
  \newtheorem{lemma}[theo]{\bfseries \hs Lemma}
  \newtheorem{corol}[theo]{\bfseries \hs Corollary}
  
  \newtheorem{algo}[theo]{\bfseries \hs Algorithm}

  \numberwithin{equation}{section} 

 \renewcommand{\span}{\mathrm{span}}

\begin{document}

 \title{Best subspace tensor approximations}
 \author{
 S. Friedland\footnotemark[1] \and
 V. Mehrmann\footnotemark[2]
}
\renewcommand{\thefootnote}{\fnsymbol{footnote}}
\footnotetext[1]{
Department of Mathematics, Statistics and Computer Science,
University of Illinois at Chicago, Chicago, Illinois 60607-7045,
USA,
\texttt{friedlan@uic.edu},
Partially supported by \emph{Berlin Mathematical School}, Berlin, Germany.
}
\footnotetext[2]{
Institut f\"ur Mathematik, TU Berlin, Str. des 17. Juni 136, D-10623 Berlin, FRG.
\texttt{mehrmann@math.tu-berlin.de}.
Partially supported by {\it Deutsche Forschungsgemeinschaft},
through the DFG Research Center {\sc Matheon}
{\it Mathematics for Key Technologies} in Berlin
}
\renewcommand{\thefootnote}{\arabic{footnote}}
\date{27 May, 2008 }
 \maketitle
 \begin{abstract}
 In many applications such as data compression, imaging or
genomic data analysis, it is important to approximate a given
tensor by a tensor that is sparsely representable. For
matrices, i.e. $2$-tensors, such a representation can be
obtained via the singular value decomposition which allows to
compute the best rank $k$ approximations. For $t$-tensors with
$t>2$ many generalizations of the singular value decomposition
have been proposed to obtain low tensor rank decompositions. In
this paper we will present a different approach which is based
on best subspace approximations, which  present an alternative
generalization of the singular value decomposition to tensors.
 \end{abstract}



 \section{Introduction}\label{intro}

In this paper we will consider data sparse approximations of
tensors. We will discuss a generalization of the singular value
decomposition from matrices to tensors that is an alternative
to the Tucker decomposition \cite{deLdV00,Tuc66}. In order not
to overload the paper with technical we will mainly discuss
$3$-tensors, but our approach will work for arbitrary tensors.

 Let $\F$ be either the field of real numbers $\R$ or complex
 numbers $\C$.  Denote by $\F^{m_1\times\ldots\times m_d}:=\otimes_{i=1}^d \F^{m_j}$
 the tensor products of $\F^{m_1},\ldots,\F^{m_d}$.
 $\cT=[t_{i_1,\ldots, i_d}]\in \F^{m_1\times\ldots\times m_d}$
 is called a \emph{$d$-tensor} in the given tensor product.  Note that the
 number of coordinates of $\cT$ is $N=m_1\ldots m_d$.
 A tensor $\cT$ is called a \emph{sparsely representable tensor} if it can
 represented with a number of coordinates that is much smaller than $N$.

 The best known example of a sparsely representable $2$-tensor
 is a low rank approximation of a matrix
 $A\in \F^{m_1\times m_2}$.  A rank $k$  approximation of $A$ is
 given by $A_{\textrm{appr}}:=\sum_{i=1}^k \u_i\v_i\trans$, which can be identified
 with $\sum_{i=1}^k \u_i\otimes \v_i$.  To store
 $A_{\textrm{appr}}$ we need only the $2k$ vectors
 $\u_1,\ldots,\u_k\in\F^{m_1},\;\v_1,\ldots,\v_k\in\F^{m_2}$.
 The best rank $k$ approximation of  $A\in \F^{m_1\times m_2}$ can be computed via
 the  \emph{singular value decomposition}, abbreviated here as SVD, \cite{GolV96}.

 The computation of the SVD requires $\mathcal {O}(m_1m_2^2+m_2^2)$
 operations and at least $\mathcal{O} (m_1 m_2)$ storage. Thus,
 if the dimensions $m_1$ and $m_2$ are very large, then the
 computation of the SVD is often infeasible. In this case other
 type of low rank approximations are considered, see e.g.
 \cite{FriKNZ06,FMMN08,GorT01}. 

 For $d$-tensors with $d>2$, however the situation is rather
 unsatisfactory. It is a major theoretical and computational
 problem to formulate good generalizations of low rank
 approximation for tensors and to give efficient
 algorithms to compute these approximations, see e.g.
 \cite{deLdV00,Lim07,Tuc66}.
 It is the goal of this paper to
 present and analyze an alternative generalization of the SVD to
 tensors. 

 A tensor $\cT=[t_{i,j,k}]\in \F^{m_1\times m_2\times m_3}$ is called a
 \emph{rank $1$} tensor, and denoted by $\cT=\u\otimes \v\otimes\w$, if
 $t_{i,j,k}=u_i v_j w_k$, where $\u=(u_1,\ldots,u_{m_1})\trans,
 \v=(v_1,\ldots,v_{m_2})\trans, \w=(w_1,\ldots,w_{m_3})\trans$.
 A tensor $\cT\in \F^{m_1\times m_2\times m_3}$ is said to have \emph{rank} $k$ if $\cT$ can
 be represented as a sum of $k$ rank $1$ tensors, and cannot be
 represented as a sum of $k-1$ rank $1$ tensors.
 Note that if $\cT$ is a sum of $k$ rank $1$ tensors, then $\cT$ can be represented
 with at most $\mathcal {O}( k(\ell+m+n))$ storage.

 We denote by $\cR(k;m_1,m_2,m_3)$ the set of tensors in
 $\F^{m_1\times m_2\times m_3}$ of rank $k$ at most.  It is easy to show that
 $\cR(1;m_1,m_2,m_3)$ is a closed set, more precisely an algebraic variety,
 in $\F^{m_1\times m_2\times m_3}$. However, it is well known,
see e.g. \cite{Fr08},
 that for some values of  $k\ge 2$, $\cR(k;m_1,m_2,m_3)$ is not a closed
 set.  (\emph{$\cR(k;m_1,m_2,m_3)$ is called a quasi-algebraic variety.})

 Let $\|\cdot\|$ be a norm on $\F^{m_1\times m_2\times m_3}$.  Then
 for $k\ge 2$ it is possible that the minimization problem
 \begin{equation}\label{brnkapp}
 \min_{\cX\in \cR(k;m_1,m_2,m_3)} \|\cT-\cX\|
 \end{equation}
 does not have a minimal solution.  This will happen if $\cT$
 has rank greater than $k$ and $\cT$ lies in the closure of
 $\cR(k;m_1,m_2,m_3)$.
  Hence, any algorithm which tries to find
 a solution to the minimization problem (\ref{brnkapp}) will fail
 for certain tensors $\cT$.  Since  $\cR(k;m_1,m_2,m_3)$
 is a closed set, for $k=1$, i.e. for the best approximation
 by a rank $1$ tensor, (\ref{brnkapp}) will always have a
 minimal solution.

 The object of this paper to introduce a new family of sparsely representable
 approximations to tensors, which we call \emph{best subspace tensor approximation (BSTA)} of a
 given tensor $\cT$.  As for the best rank $1$ approximation, we will show that the BSTA always exists.
 Due to this  fact, we think that in the case
 that the norm $\|\cdot\|$ on $\F^{m_1\times m_2\times m_3}$ is the norm induced by the inner products
 on the vector spaces $\F^{m_1},\;\F^{m_2},\;\F^{m_3}$,  the BSTA is
 an appropriate generalization of the SVD,
 see \cite{deLdV00} for other generalizations  of the SVD for tensors.
 Similar approach was suggested recently by
 Khoromskij \cite{Kho}.
 We will also present a numerical algorithm
to compute the best subspace tensor approximation
 that is based on the computation of singular value decompositions for matrices.

 Unfortunately this numerical algorithm is extremely expensive.
 In order to reduce the complexity, in the last section we
 consider a procedure that is based on the recently suggested
 fast SVD \cite{FMMN08}.

 \section{Notation and preliminary results}\label{notprem}

 We denote by a bold capital letter a finite dimensional vector
 space $\U$ over the field $\F$.  A vector $\u\in\U$ is denoted
 by a bold face lower case letter.  A matrix
 $A\in \F^{m_1\times m_2}$ denoted by a capital letter $A$, and we let either
 $A=[a_{i,j}]_{i=j=1}^{m_1\times m_2}$ or simply $A=[a_{i,j}]$.
 A \emph{3-tensor array}
 ${\cT}\in \F^{m_1\times m_2\times m_3}$ will be
 denoted by a capital calligraphic letter.  So either
 ${\cT}=[t_{i,j,k}]_{i=j=k=1}^{m_1,m_2,m_3}$ or simply ${\cT}=[t_{i,j,k}]$.
For a positive integer $n$ we also use the convenient notation $\an{ n}:=\{1,2,\ldots,n\}$.

 Let $\U_1,\U_2,\U_3$ be three vectors spaces over $\F$ with $m_j:
 =\dim \U_j$, $j=1,2,3$ and let $\u_{1,j},
 \ldots,\u_{m_j,j}$ be a basis of $\U_j$ for $j=1,2,3$.
 Then $\U:=\U_1\otimes\U_2\otimes\U_3$ is the \emph{tensor product}
 of $\U_1$, $\U_2$, and $\U_3$; $\U$ is a vector space of dimension
 $m_1m_2m_3$, and
 \begin{equation}\label{tenbas}
 \u_{i_1,1}\otimes\u_{i_2,2}\otimes\u_{i_3,3},\quad
 i_j= 1,\ldots,m_j,\ j=1,2,3,
 \end{equation}
 is a basis of $\U$.

 A \emph{$3$-tensor} $\tau$ is a vector in $\U$ and it has a representation
 \begin{equation}\label{tenrep}
 \tau=\sum_{i_1=i_2=i_3=1}^{m_1,m_2,m_3}
 t_{i_1,i_2,i_3}\u_{i_1,1}\otimes\u_{i_2,2}\otimes \u_{i_3,3},
 \end{equation}
 in the basis (\ref{tenbas}).  If the basis (\ref{tenbas}) is fixed
 then $\tau$ is identified with ${\cT}=[t_{i_1,i_2,i_3}]\in
 \F^{m_1\times m_2\times m_3}$.

 Recall that $\x_1\otimes\x_2\otimes \x_3$, were $\x_i\in\U_i,
 i=1,2,3$, is called a \emph{rank $1$} tensor.
 (Usually one assumes that all $\x_i\ne \0$.  Otherwise $\0=\x_1\otimes\x_2\otimes\x_3$ is called
 a rank $0$ tensor.) Then
 (\ref{tenrep}) is a decomposition
 of $\tau$ as a sum of at most $m_1m_2m_3$  rank $1$ tensors, as
 $t_{i_1,i_2,i_3}\u_{i_1,1}\otimes\u_{i_2,2}\otimes \u_{i_3,3}=
 (t_{i_1,i_2,i_3}\u_{i_1,1})\otimes\u_{i_2,2}\otimes \u_{i_3,3}$.
 A decomposition of $\tau\in \U \backslash \{\0\}$ as a sum of rank
 $1$ tensors is given  by
 \begin{equation}\label{rankonedec}
 \tau=\sum_{i=1}^k \x_i\otimes\y_i\otimes \z_i,\quad
 \x_i\in\U_1,\; \y_i\in\U_2,\; \z_i\in\U_3, \;i=1,\ldots,k.
 \end{equation}
 The minimal $k$ for which the above equality holds is called
 the \emph{rank} of the tensor $\tau$.  This definition is completely analogous
 to the definition of the rank for a matrix $A=[a_{i_1,i_2}]\in \F^{m_1\times m_2}$, which can be
 identified with 2-tensor in $\sum_{i_1=i_2=1}^{m_1,m_2} a_{i_1,i_2} \u_{i_1,1}\otimes
 \u_{i_2,2}\in\U_1\otimes \U_2$.

 For $j\in\{1,2,3\}$  denote by
 $j^c:=\{p,q\}=\{1,2,3\}\backslash\{j\}$, where $1\le p <q\le 3$,
and set $\U_{j^c}=\U_{\{p,q\}}:=\U_p\otimes\U_q$.

A tensor $\tau\in \U_1\otimes\U_2\otimes \U_3$ induces a linear
 transformation $\tau(j):\U_{j^c}\to \U_j$ as follows.  Suppose
 that $\u_{1,\ell},\ldots,\u_{m_\ell,\ell}$ is a basis in $\U_\ell$ for
 $\ell=1,2,3$.  Then any $\v\in \U_{j^c}$ is of the form
 $$\v=\sum_{i_p=i_q=1}^{m_p,m_q} v_{i_p,i_q} \u_{i_p,p}\otimes
 \u_{i_q,q}$$
and the application of $\tau(j)$ is given by
 \begin{equation}\label{deftauj}
 \tau(j)\;\v=\sum_{i_j=1}^{m_j}\big(\sum_{i_p,i_q=1}^{m_p,m_q} t_{i_1,i_2,i_3}
 v_{i_p,i_q}\big)\u_{i_j,j}.
 \end{equation}
 Then $\rank_j (\tau)$ \emph{is the rank of the operator} $\tau(j)$.
 Equivalently, let $A(j)=[a_{\ell,i_j}]\in \R^{m_pm_q\times m_j}$, where
 each integer $\ell\in \an{ m_pm_q}$ corresponds to a pair $(i_p,i_q)$,
 for $ i_p=1,\ldots,m_p,\ i_q=1,\ldots,m_q$, and $i_j\in \langle m_j\rangle$.
 (For example we may arrange
 the pairs $(i_p,i_q)$ in the lexicographical order.
 Then $i_p=\lceil \frac{\ell}{m_q}\rceil$
 and $i_q=\ell-(i_p-1)m_q$.)  Set $a_{\ell,i_j}=t_{i_1,i_2,i_3}$.  Then
 $\rank_j(\tau)=\rank\; A(j)$.

 The following proposition is straightforward.
 \begin{prop}\label{rankdef}
 Let $\tau\in \U_1\otimes\U_2\otimes\U_3$ be given by
 (\ref{tenrep}).  Fix $j\in\{1,2,3\}$ and set $ j^c=\{p,q\}$.  Let
 $T_{i_j,j}:=[t_{i_1,i_2,i_3}]_{i_p=i_q=1}^{m_p,m_q}\in \F^{m_p\times
 m_q}, i_j=1,\ldots,m_j$.  Then $\rank_j(\tau)$ is the dimension of subspace
 of $m_p\times m_q$ matrices spanned by $T_{1,j},\ldots,T_{m_j,j}$.
 \end{prop}

 Assume that each $\U_j$ is an inner
 product space, with the inner product $\an{\cdot,\cdot}_j$
 for $j=1,2,3$.  Let $\u_{1,j},\ldots,\u_{m_j,j}$,  $j=1,2,3$ be an
 orthonormal basis in $\U_j$ with respect to $\an{\cdot,\cdot}_j$.  Define
 an inner product on $\U$, denoted by $\an{\cdot,\cdot}$,  by assuming that the basis
 (\ref{tenbas})  is an orthonormal basis in $\U$.  It is straightforward to
 show that the above inner product does not depend on the
 choice of the orthonormal bases in $\U_1,\U_2,\U_3$.
 The so defined  inner product in $\U$ is called the \emph{induced
 inner product} and we have identity
 $$\an{\x\otimes\y\otimes \z,\u\otimes\v\otimes\w}=\an{\x,\u}_1
 \an{\y,\v}_2\an{\z,\w}_3.$$

 On $\F^{m_1\times m_2\times m_3}$ the standard inner product
 $\an{\cX,\cY}$ is given by $\sum_{i=j=k}^{m_1,m_2,m_3} x_{i,j,k}\bar
 y_{i,j,k}$, where $\cX=[x_{i,j,k}],\cY=[y_{i,j,k}]$.  This inner product is induced
 by the standard inner products on $\F^{m_1},\F^{m_2},\F^{m_3}$.  So
 $\|\cX\|=(\sum_{i=j=k=1}^{m_1,m_2,m_3} |x_{i,j,k}|^2)^{\frac{1}{2}}$
 is the \emph{Hilbert-Schmidt norm} on $\F^{m_1\times m_2\times m_3}$.

 We denote by $\Gr(p,\F^n)$ the set of all $p$-dimensional
 subspaces of $\F^n$.  It is well known that $\Gr(p,\F^n)$ is
 a closed set, more precisely an algebraic variety, called
 the \emph{Grassmannian} of $\F^n$ \cite{Har92}.

 \begin{defn}\label{deftpg}
 Let $p\in \an{m_1}, q\in \an{m_2}, r\in \an{m_3}$.  Denote by
 $\Gr(p,\F^{m_1})\otimes\Gr(q,\F^{m_2})\otimes\Gr(r,\F^{m_3})\subseteq
 \Gr(pqr,\F^{m_1\times m_2\times m_3})$ the set of all
 $pqr$-dimensional subspaces in $\F^{m_1\times m_2\times m_3}$
 of the form $\X\otimes\Y\otimes\Z$, where
 $\X\in\Gr(p,\F^{m_1}),\Y\in\Gr(q,\F^{m_2}),\Z\in\Gr(r,\F^{m_3})$.
 \end{defn}

 Clearly, $\Gr(p,\F^{m_1})\otimes\Gr(q,\F^{m_2})\otimes\Gr(r,\F^{m_3})$
 is a closed subvariety of
 \par\noindent
 $\Gr(pqr,\F^{m_1\times m_2\times m_3})$.
Define by  $\dist(\cT,\rS, \|\;\|):=\inf_{\cX\in\rS} \|\cT-\cX\|$
the distance of $\cT$ to a set $\rS\subset \F^{m_1\times m_2 \times m_3}$
 with respect to the norm $\|\;\|$.  Then the \emph{best $(p,q,r)$
 subspace approximation of $\cT\in \F^{l\times m \times n}$} is
 given by
 \begin{equation}\label{bpqrsap}
 \min_{\X\otimes\Y\otimes\Z\in
 \Gr(p,\F^{m_1})\otimes\Gr(q,\F^{m_2})\otimes\Gr(r,\F^{m_3})}
 \dist(\cT,\X\otimes\Y\otimes\Z,\|\;\|),
 \end{equation}
and we denote the subspace where the minimum is achieved by $
\X^*\otimes\Y^*\otimes\Z^*$ and the minimal tensor by $\cX^*\in
 \X^*\otimes\Y^*\otimes\Z^*$, i.e. we have
\begin{equation}\label{minargtens}
 \dist(\cT,\X^*\otimes\Y^*\otimes\Z^*,\|\;\|)=\|\cT-\cX^*\|.
\end{equation}
%

 Let $\ell_1\in \an{m_1}, \ell_2\in \an{m_2}, \ell_3\in \an{m_3}$ and suppose that $\U_1\in
 \Gr(\ell_1,\F^{m_1}),\U_2\in\Gr(\ell_2,\F^{m_2}),\U_3\in\Gr(\ell_3,\F^{m_3})$.
 Choose
 $$
 \u_{1,1},\ldots,\u_{\ell_1,1}\in\F^{m_1},\;\u_{1,2},\ldots,\u_{\ell_2,2}\in
 \F^{m_2},\;\u_{1,3},\ldots,\u_{\ell_3,3}
 \in\F^{m_3},
 $$
 such that $\u_{1,j},\ldots,\u_{\ell_j,j}$ is an orthonormal basis
 in $\U_j$ for $j=1,2,3$.  Then for $\tau\in \F^{m_1\times m_2\times m_3}$
 let
 \begin{equation}\label{cordt}
 t_{i,j,k}=\an{\tau,\u_{i,1}\otimes\u_{j,2}\otimes\u_{k,3}}, \quad
 i=1,\ldots,l,\;j=1,\ldots,m,\;k=1,\ldots,n.
 \end{equation}
 So $\cT=[t_{i,j,k}]$ is the representation of $\tau$ in the
 orthonormal basis.
 Then
 \begin{equation}\label{projt}
 P_{\U_1\otimes\U_2\otimes\U_3}
 (\tau)=\xi=\sum_{(i,j,k)\in\an{\ell_1}\times\an{\ell_2},\times\an{\ell_3}}
 t_{i,j,k} \u_{i,1}\otimes\u_{j,2}\otimes\u_{k,3}
 \end{equation}
 is the orthogonal projection of $\tau$ on the subspace
 $\U_1\otimes\U_2\otimes\U_3$.
 Thus
 \begin{equation}\label{dishs}
 \dist(\tau,\U_1\otimes\U_2\otimes\U_3)=\|\tau-\xi\|=(\sum_{(i,j,k)\in
 \an{m_1}\times\an{m_2}\times\an{m_3}\backslash
 \an{\ell_1}\times\an{\ell_2}\times\an{\ell_3}} |t_{i,j,k}|^2)^{\frac{1}{2}}
 \end{equation}
is  the distance with respect to the Hilbert-Schmidt norm on $\F^{m_1\times m_2\times m_3}$.
 Clearly, we have
 \begin{equation}\label{normeq}
 \|\tau\|^2=\|P_{\U_1\otimes\U_2\otimes\U_3}(\tau)\|^2
 +\dist(\tau,\U_1\otimes\U_2\otimes\U_3)^2.
 \end{equation}

 \section{The SVD as best subspace tensor approximation}\label{svdbsta}
In this section we will illustrate that the SVD allows to compute
the best subspace tensor approximation for $2$-tensors.

 Let us view $m_1\times m_2$ matrices  as
 $2$-tensors.  Here $\x\otimes \y$ corresponds to the matrix
 $\x\y\trans$.  A tensor $\tau\in \F^{m_1}\otimes\F^{m_2}$ can be viewed as
 a linear transformation $\tau:\F^{m_1}\to \F^{m_2}$ as follows.
 First observe that a rank $1$ tensor $\x\otimes \y$ gives rise to the
 linear transformation $(\x\otimes\y)(\z)=\an{\z,\bar\y}\x$.
 Now extend this notion to any $\tau\in\F^{m_1}\otimes\F^{m_2}$, which
 is a sum of rank $1$ tensors.

 We claim that the best rank $k$ approximation of $\tau$ is obtained
as the solution  to  the minimization problem
 \begin{equation}\label{matsvd}
 \min_{\X\in \Gr(k,\F^{m_1}),\Y\in\Gr(k,\F^{m_2})}
 \dist(\tau,\X\otimes\Y)=\dist(\tau,\X^*\otimes\Y^*),
 \end{equation}
where $\X^*,\Y^*$ are the subspaces spanned by the $k$ left and
 right singular vectors of $\tau$ associated with the largest $k$ singular values.

 Indeed, suppose that the minimum in (\ref{matsvd}) is achieved
 for some tensor $\alpha\in \X^*\otimes \Y^*$, so $\rank\; \alpha \le k$.
Hence the best approximation by a rank
 $k$ tensor is not worse than the minimum of (\ref{matsvd}).  On the other
 hand, any rank $k$ tensor is an element of sum $\X\otimes\Y$ for
 some $\X\in\Gr(k,\F^{m_1}),\Y\in\Gr(k,\F^{m_2})$.  So the minimum in
 (\ref{matsvd}) is not bigger than the best rank $k$
 approximation.  But the best rank $k$ approximation to a given $2$-tensor is obtained by the
 SVD \cite{GolV96}.

 We now consider the following approximation problems for $2$-tensors, which
 is equivalent to the corresponding matrix problem.
%
 \begin{lemma}\label{minprobmat}  Let $\Y\subset \F^{m_2}$ be a
 given $\ell_1\in \an{m_1}$ dimensional subspace.  For $i\in\an{m_1}$
 and $\tau\in \F^{m_1}\otimes \F^{m_2}$ consider the minimization problem
of finding $X\in \Gr(i,\F^{m_1})$ such that
\begin{equation}\label{matsvdm}
 \min_{\X\in \Gr(i,\F^{m_1})}
 \dist(\tau,\X\otimes\Y)=\dist(\tau,\X^*\otimes\Y).
 \end{equation}
View $\tau$ as a linear mapping from $\F^{m_1}$ to $\F^{m_2}$.
If $\dim (\tau \Y)\le i$ then $\X^*$ is any subspace that contains
 $\tau \Y$.  If $\dim(\tau \Y)>i$ then $\X^*$ is the subspace
 spanned by the left singular vectors associated with
the $i$ largest singular values of $\tau|_\Y$ (which is a linear map $\tau:\Y\to \F^{m_1}$).
 \end{lemma}
 \proof  Choose the standard orthonormal basis
 $\e_1,\ldots,\e_{m_1}\in\ F^{m_1}$ and an orthonormal basis
 $\y_1,\ldots,\y_{m_2}\in\F^{m_2}$ such that $\Y=\span
 (\y_1,\ldots,\y_\ell)$.  Let
 $\Y^{\perp}=\span(\y_{\ell+1},\ldots,\y_{m_2})$.
 Then $\F^{m_1}\otimes\F^{m_2}=\F^{m_1}\otimes \Y \oplus
 \F^{m_1}\otimes\Y^{\perp}$ is an orthogonal decomposition of
 $\F^{m_1}\otimes\F^{m_2}$. This means that we can write $\tau$ as
 $$
\tau=\phi+\psi,\;
 \phi=P_{\F^{m_1}\otimes\Y}(\tau),\;\psi=P_{\F^{m_1}\otimes\Y^{\perp}}(\tau),\;
 \|\tau\|^2=\|\phi\|^2+\|\psi\|^2.
$$
 Since we require $\X\otimes\Y\subset \F^{m_1}\otimes \Y$ it follows that the
 minimization problem (\ref{matsvdm}) is equivalent to the minimization
 problem
 \begin{equation}\label{matsvdm1}
 \min_{\X\in \Gr(i,\F^{m_1})}
 \dist(\phi,\X\otimes\Y)=\dist(\tau,\X^*\otimes\Y).
 \end{equation}
 Observe next that $\phi$, viewed as a linear transformation
 $\phi:\F^{m_1}\to \Y$ is equal to $\tau|\Y$.
 The classical result for matrices implies that the best rank ´$i$  approximation of $\phi$ is given via
 the left singular vectors associated to the largest $i$ singular values of $\phi$.
 \qed

 In this section we have shown that the best subspace tensor approximation for $2$-tensors
 is obtained via the singular value decomposition. This immediately suggest to use
 it as a generalization of the SVD for higher tensors.

 \section{Best subspace tensor approximations for $3$-tensors}

 I n this section we study the best subspace tensor
 approximation for $3$-tensors.
 Let $\tau\in \F^{m_1\times m_2\times m_3}$ and assume that $p\in
 \an{m_1},q\in \an{m_2},r\in \an{m_3}$ and consider the minimization problem
 \begin{equation}\label{bstahs}
 \min_{\X\otimes\Y\otimes\Z\in
 \Gr(p,\F^{m_1})\otimes\Gr(q,\F^{m_2})\otimes\Gr(r,\F^{m_3})} \dist(\tau,\X\otimes\Y\otimes\Z)
\end{equation}
 and suppose that is minimum is achieved for the subspace
 $\X^*\otimes\Y^*\otimes\Z^*$ with the tensor
 $\xi$, i.e.
\[
 \dist(\tau,\X^*\otimes\Y^*\otimes\Z^*)=\|\tau-\xi^*\|,\;\xi^*\in
 \X^*\otimes\Y^*\otimes\Z^*.
 \]
 In view of (\ref{normeq}) this minimization problem is
equivalent to the maximization problem
  \begin{equation}\label{bstahs1}
 \max_{\X\otimes\Y\otimes\Z\in
 \Gr(p,\F^{m_1})\otimes\Gr(q,\F^{m_2})\otimes\Gr(r,\F^{m_3})}
 \|P_{\X\otimes\Y\otimes\Z}(\tau)\|^2=\|
 P_{\X^*\otimes\Y^*\otimes\Z^*}(\tau)\|^2.
 \end{equation}
 To simplify our exposition we state our results for
 $\F=\R,\C$, but we give the proofs only for $\F=\R$.

 To solve the minimization problem, we study the critical points (i.e. the
 points of vanishing gradient) of
 $\|P_{\X\otimes\Y\otimes\Z}(\tau)\|^2$ on
 $\Gr(p,\F^{m_1})\otimes\Gr(q,\F^{m_2})\otimes\Gr(r,\F^{m_3})$.
 To do that we need the following lemma which follows from the
 Courant-Fischer theorem, see e.g. \cite{GolV96}.
 In the following, we use  $\Fr(i,\F^{m_1})$ to denote the
 manifold of all sets of  $i$ orthonormal vectors  $\{\x_1,\ldots,\x_i\}\subset\F^{m_1}$.
 \begin{lemma}\label{critptkfn}
 Let $B\in\F^{m_1\times m_1}$
 be a Hermitian matrix.  Let  a linear functional $g_B:\Fr(i,\R^{m_1})\to \R$ be given
 by $g_B(\x_1,\ldots,\x_i)=\sum_{l=1}^i \x_l\trans B\x_l$.
 Then the critical points of $g_B$ are all sets
 $\{\x_1,\ldots,\x_i\}$ such that $\span(\x_1,\ldots,\x_i)$
 contains $i$ linearly independent eigenvectors of $B$.
 \end{lemma}

 \proof
 We prove the lemma by induction on $i$.  For $i=1$ we have
 $g_B(x)=\x\trans B\x$ (note that $\| \x\|=1$).
 Then by the Courant-Fischer Min-Max characterization, see e.g. \cite{GolV96},
 $\x\ne\0$ is a critical point if and only if $\x$
 is an eigenvector of $B$.

 If $\x_1,\ldots,\x_i$ are eigenvectors of $B$ it is
 straightforward to see that $\{x_1,\ldots,\x_i\}$ is a
 critical point of $g_B$.  Indeed, consider a variation
 $x_\ell(t)=\x_\ell+t\u_\ell+t\v_\ell +O(t^2), \ell=1,\ldots,i$, where
 $\u_l\in \span(\x_1,\ldots,\x_i),\ \v_l\in
 \span(\x_1,\ldots,\x_i)^{\perp}$.  Then the contribution
 involving $\u_1,\ldots,\u_i$ is quadratic in $t$.
 Since $\v_\ell\trans \x_l=0,\ \ell=1,\ldots,i$ it follows that the
 contribution in $\v_1,\ldots,\v_i$ is also quadratic in $t$.
 It remains to show that if $\{\x_1,\ldots,\x_i\}$ is a
 critical point of $g_B$ then $\span(\x_1,\ldots,\x_i)$ is
 spanned by $i$ eigenvectors of $B$.

 Suppose that the assertion holds for $i=k-1$ and
 assume that $i=k\le m_1$.  If $k=m_1$ then the assertion is
 clear because the whole space is spanned
by eigenvectors of $B$.  So let $k<m$.
 Note that if
 $\{\y_1,\dots,\y_i\}\in\Fr(i,\R^{m_1})$ and
 $\span(\y_1,\ldots,\y_i)=\span(\x_1,\ldots,\x_i)$ then
 $g_B(\x_1,\ldots,\x_i)=g_B(\y_1,\ldots,\y_i)$.
 So we may  assume w.l.o.g. that the matrix
\[
C=[\x_s\trans
 B\x_t]_{s,t=1}^i\]
 is diagonal. Furthermore, we may assume that
 $\x_s=\e_s, s=1,\ldots,i$.  The induction hypothesis states
 that for any $k\in \{i+1,\ldots,m_1\}$ the symmetric matrix $B_k$,
 obtained by erasing $k$ rows and columns of $B$ is a direct sum
 of $C$ and the corresponding other block.  Hence $B=C\oplus
 C'$ and the assertion follows.
 \qed

We immediately have the following corollary.
%
 \begin{corol}\label{crit2tenapr}  Let $\alpha\in\F^{m_1}\otimes
 \F^{m_2}$ and suppose that $i$ is an integer in the interval $[1,m_1]$.
 Then $\U\in\Gr(i,\F^{m_1})$ is a critical point of the linear functional
 $\|P_{\X\otimes\F^{m_2}}(\alpha)\|^2: \Gr(i,\F^{m_1})\to [0,\infty)$
 if and only if $\U$ is spanned by some $i$ left singular vectors
 of the induced dual operator $\tilde \alpha: \F^{m_2}\to\F^{m_1}$.
 (Here some singular vectors may correspond to the singular value $0$.)
 \end{corol}
 \proof Represent $\tilde \alpha$ by $A\in \R^{m_1\times m_2}$
and let $B=A A\trans$. Let $\X\in\Gr(i,\R^{m_1})$ and suppose that
 $\{\x_1,\ldots,\x_i\}\in\Fr(i,\R^m)$ is a basis of $\X$.
 Then $\|P_{\X\otimes\F^{m_2}}(\alpha)\|^2=g_B(\x_1,\ldots,\x_i)$, and the result
follows from Lemma~\ref{critptkfn}.
\qed

We will now construct projections of $3$-tensors to $2$-tensors, which we can use to
compute best subspace approximations.

 Let $\tau\in \F^{m_1}\otimes\F^{m_2}\otimes\F^{m_3}$ and
 $\X\in\Gr(p,\F^{m_1}),\Y\in\Gr(q,\F^{m_2}),\Z\in\Gr(r,\F^{m_3})$.
Suppose that
 $\e_1,\ldots,\e_{m_1},\;\f_1,\ldots,\f_{m_2},\;\g_1,\ldots, \g_{m_3}$
 are orthonormal bases in $\F^{m_1},\F^{m_2},\F^{m_3}$ respectively, such
 that $\e_1,\ldots,\e_p,\;\f_1,\ldots,\f_q,\;\g_1,\ldots, \g_r$
 are bases of $\X,\Y,\Z$, respectively.
 Then we can express $\tau$ as $\tau=\sum_{i=j=k=1}^{m_1,m_2,m_3} t_{i,j,k}
 \e_i\otimes\f_j\otimes\g_k$ and consider the following linear operators.
 \begin{enumerate}
 \item
The first operator
 $\tau(\Y,\Z):\F^{m_1}\to \Y\otimes \Z$ is constructed
 as follows.  View $P_{\F^{m_1}\otimes\Y\otimes\Z}(\tau)$ as a tensor
 in $\F^{m_1}\otimes\Y\otimes\Z$, i.e.
 $$P_{\F^{m_1}\otimes\Y\otimes\Z}(\tau)
 =\sum_{i=j=k=1}^{m_1,q,r} t_{i, j, k}
 \e_{i}\otimes\f_{j}\otimes\g_{k}$$
and then define for $\x\in\F^{m_1}$ the operator via
 $$\tau(\Y,\Z)(\x)=\sum_{i=j=k=1}^{m_1,q,r} t_{i,j,k}\an{\x,\e_{i}}_1\
 \f_{j}\otimes\g_{k},$$
where as before $\an{\cdot, \cdot}_1$ denotes the inner product in $\F^{m_1}$.
 \item Analogously we proceed for $\tau(\X,\Z):\F^{m_2}\to \X\otimes \Z$.
We view $P_{\X\otimes\F^{m_2}\otimes\Z}(\tau)$ as a tensor
 in $\X\otimes\F^{m_2}\otimes\Z$, i. e.,
 $$P_{\X\otimes\F^{m_2}\otimes\Z}(\tau)
 =\sum_{i=j=k=1}^{p,m_2,r} t_{i,j,k}
 \e_{i}\otimes\f_{j}\otimes\g_{k}$$
and then for any $\y\in\F^{m_2}$ we define the operator via
 $$
\tau(\X,\Z)(\y)=\sum_{i=j=k=1}^{p,m_2,r} t_{i,j,k}\an{\y,\f_{j}}_2\
 \e_{i}\otimes\g_{k}.$$
 \item Finally $\tau(\X,\Y):\F^n\to \X\otimes \Y$
 is given as follows.  View $P_{\X\otimes\Y\otimes\F^{m_3}}(\tau)$ as a tensor
 in $\X\otimes\Y\otimes\F^{m_3}$, i. e.,
 $$P_{\X\otimes\Y\otimes\F^{m_3}}(\tau)
 =\sum_{i=j=k=1}^{p,q,m_3} t_{i,j,k}
 \e_{i}\otimes\f_{j}\otimes\g_{k}.$$
 Then for any $\z\in\F^{m_3}$, we define the operator via
 $$\tau(\X,\Y)(\z)=\sum_{i=j=k=1}^{p,q,m_3} t_{i,j,k}\an{\z,\g_{k}}_3\
 \e_{i}\otimes\f_{j}.$$
  \end{enumerate}
We have the following theorem.
 \begin{theo}\label{charcritpfXYZ}  Let $0\ne \tau\in
 \F^{m_1}\otimes\F^{m_2}\otimes\F^{m_3}$.  Let $p\in\an{m_1},q\in\an{m_2},r\in\an{m_3}$.
 Then $\U\in \Gr(p,\F^{m_1}),
 \V\in\Gr(q,\F^{m_2}),\W\in\Gr(r,\F^{m_3})$ is a critical point of
 $\|P_{\X\otimes\Y\otimes\Z}(\tau)\|^2$ on
 $\Gr(p,\F^{m_1})\otimes\Gr(q,\R^{m_2})\otimes\Gr(k,\F^{m_3})$ if and only if
 the following conditions hold
 \begin{enumerate}
 \item $\U$ is spanned by some $p$ left singular vectors of
 $\tau(\V,\W)$.
 \item $\V$ is spanned by some $q$ left singular vectors of
 $\tau(\U,\W)$.
 \item $\W$ is spanned by some $r$ left singular vectors of
 $\tau(\U,\V)$.
 \end{enumerate}
 \end{theo}
 \proof
 Since the critical points are the zeros of the first
 derivative, it is enough to prove the necessary
 conditions for the function $\|P_{\X\otimes\V\otimes\W}(\tau)\|^2$.
 Considering this as  function on $\Gr(p,\R^{m_1})$, Condition
 1. then  follows immediately by Corollary~\ref{crit2tenapr}.
 The other conditions follow analogously.
 \qed

 In the following we will describe an iterative procedure to compute the best subspace tensor
 approximation. In order to find good starting values for $\U=\X_0,\V=\Y_0,\W=\Z_0$ we make use
 of the  SVD.  As explained in \S2 we can \emph{unfold} $\tau$ as a
 matrix $A_1$, say $m_1\times (m_2n_3)$, by considering $\tau(1)$
 as defined in (\ref{deftauj}).  Then we perform the SVD
and use as approximation the corresponding $p$-dimensional $\X_0\in
 \Gr(p,\F^{m_1})$ spanned the left singular vectors of
 $A_1$ associated with the $p$ largest singular values.  In a
similar way we determine $\Y_0\in\Gr(q,\F^{m_2}),\Z_0\in
 \Gr(r,\F^{m_3})$.

To find the maximum in
 (\ref{bstahs1}) we then apply a relaxation method.
 \begin{algo}\label{relaxalgo3ten} Let
 $\tau\in\F^{m_1}\otimes\F^{m_2}\otimes\F^{m_3}$, $p\in\an{m_1},q\in\an{m_2},r\in\an{m_3}$ and
starting values $\X_0\in\Gr(p,\F^{m_1}),\Y_0\in\Gr(q,\in\F^{m_2}),\Z_0\in\Gr(r, \F^{m_3})$ be given.

Suppose  that $(\X_i,\Y_i,\Z_i)$ have been computed.  Then
 \begin{enumerate}
 \item $\X_{i+1}$ is obtained as the $p$-dimensional subspace
 corresponding to left singular
 vectors of $\tau(\Y_i,\Z_i)$ associated with the $p$ largest singular values.
 \item $\Y_{i+1}$ is obtained as the $q$-dimensional subspace
 corresponding to the  left singular
 vectors of $\tau(\X_{i+1},\Z_i)$ associated with the $q$ largest singular values.
 \item $\Z_{i+1}$ is obtained as the $r$-dimensional subspace
 corresponding to the left singular
 vectors of $\tau(\X_{i+1},\Y_{i+1})$ associated with the $r$ largest singular values.
 \end{enumerate}
 \end{algo}

 We have the following convergence result.
 \begin{corol}\label{critptrelalg} The subspaces
 $\X_i,\Y_i,\Z_i, i=0,1,\ldots$ defined in
 Algorithm~\ref{relaxalgo3ten} converge to subspaces $\U,\V,\W$
 which give a critical point of $\|P_{\X\otimes\Y\otimes\Z}(\tau)\|^2$.
 Moreover, this critical point is a maximal point, with
 respect to any one variable, when the other variables are fixed.
 Furthermore the following conditions hold.
 \begin{enumerate}
 \item $\U$ is spanned by the left singular vectors of
 $\tau(\V,\W)$ associated with the $p$ largest values.
 \item $\V$ is spanned by the left singular vectors of
 $\tau(\U,\W)$ associated with the $q$ largest values.
 \item $\W$ is spanned by the  left singular vectors of
 $\tau(\U,\V)$ associated with the $r$ largest singular values.
 \end{enumerate}
 \end{corol}

 In this section we have shown that the best subspace tensor approximation for $3$-tensors is a
 a generalization of the singular value decomposition. It is obvious how this procedure
 can be extended to arbitrary $k$ tensors.

Unfortunately the described procedure is extremely expensive,
since in every step a singular value decomposition of a very
large full matrix has to be performed. In order to reduce the
complexity, in the next section we consider a procedure that is
based on the recently suggested fast SVD \cite{FMMN08}.

\section{Fast low rank $3$-tensors approximations}
 In this section we generalize the algorithm outlined in
 \cite{FMMN08} to the fast low rank tensor
 approximation, abbreviated as FLRTA, to $3$-tensors.
 Let $\cA=[a_{i_1,i_2,i_3}]\in
 \R^{l_1\times l_2\times l_3}$ be a $3$-tensor, where
 the dimensions $l_1,l_2,l_3$, are large.  For each $j=1,2,3$ we
 read subtensors of $\cA$ denoted by $\cC_j=[c^{(j)}_{i_{1,j} i_{2,j}
 i_{3,j}}]\in \R^{l_{1,j}\times l_{2,j} \times l_{3,j}}$.
 We assume that $\cC_j$ has the same number of coordinates
 as $\cA$ in $j$-th direction, and a small number of
 coordinates in the other two directions.
 That is, $l_{j,j}=l_j$ and the other two indices
 $l_{s,j},\ s\in \{1,2,3\}\backslash\{j\}$ are of order $\mathcal O(k)$,
 for $j=1,2,3$.  So $\cC_j$ corresponds to the $j$-\emph{section} of the tensor
 $\cA$.  The \emph{small} dimensions of $\cC_j$ are
 $(l_{s_{j},j},l_{t_{j},j})$ where $\{s_{j},t_{j}\}=\{1,2,3\}\backslash\{j\}$
 for $j=1,2,3$.  Let
 $m_j:=l_{s_{j},j}l_{t_{j},j}$ for $j=1,2,3$.

To determine an approximation, we then look for a $6$-tensor
 $$\cV=[v_{q_1, q_2, q_3, q_4, q_5, q_6}]
\in \R^{l_{2,1}\times l_{3,1}\times l_{1,2}\times l_{3,2}
 \times l_{1,3}\times l_{2,3}}$$
and approximate the given tensor $\cA$ by a tensor
 $$\cB=[b_{i_1,i_2,i_3}]:=\cV\cdot\cC_1 \cdot\cC_2\cdot\cC_3 \in
 \R^{\ell_1\times  \ell_2\times \ell_3},$$
 where we contract the $6$
 indices in $\cV$ and the corresponding two indices
 $\{1,2,3\}\backslash\{j\}$ in $\cC_j$ for $j=1,2,3$, i.e.,
our approximation has the entries
 \begin{equation}\label{tencon}
 b_{i_1,i_2,i_3}=\sum_{q_1=1}^{\ell_{2,1}}
\sum_{q_2=1}^{\ell_{3,1}}\sum_{q_3=1}^{\ell_{1,2}}\sum_{q_4=1}^{\ell_{3,2}}
\sum_{q_5=1}^{\ell_{1,3}}\sum_{q_6=1}^{\ell_{2,3}}
 v_{q_1, q_2, q_3, q_4, q_5, q_6}
 c^{(1)}_{i_1,q_1,q_2}c^{(2)}_{q_3,i_2,q_4}c^{(3)}_{q_5,q_6,i_3}.
 \end{equation}
 This approximation is  equivalent to a so-called \emph{Tucker
 approximation} \cite{Tuc66}.
Indeed, if we represent each tensor $\cC_j$ by a matrix
 $C_j\in \R^{m_j\times l_j}$ that has
 the same number of columns as the range of the $j$-th index of
 the tensor $\cA$ and as number of rows  the product of
 the ranges of the remaining two \emph{small} indices of
 $\cC_j$, i.e.  $C_j=[c^{(j)}_{r,i_j}]_{r,i_j=1}^{m_j\cdot\ell_j}$.
 Then
 $c^{(j)}_{r,i_j}$ is equal to the
 corresponding entry $c^{(j)}_{i_1,i_2,i_3}$, where the value of $r$
 corresponds to the double index $(i_s,i_t)$ for
 $\{s,t\}=\{1,2,3\}\backslash\{j\}$.

 Now with $\cU=[u_{j_1,j_2,j_3}]\in \R^{m_1\times m_2\times
 m_3}$, the equivalent Tucker representation of
 $\cB=[b_{i_1,i_2,i_3}]$ is given by the entries
 \begin{equation}\label{tucrep1}
 b_{i_1,i_2,i_3}=\sum_{j_1=1}^{m_1} \sum_{j_2=1}^{m_2}
 \sum_{j_3=1}^{m_3}
 u_{j_1,j_2,j_3} c^{(1)}_{j_1,i_1}
 c^{(2)}_{j_2,i_2}c^{(3)}_{j_3,i_3}, \quad (i_1,i_2,i_3)\in
 \an{\ell_1}\times\an{\ell_2}\times\an{\ell_3}.
 \end{equation}
 This formula is expressed commonly as
 \begin{equation}\label{tucrep2}
 \cB=\cU\times_1 C_1 \times _2 C_2 \times _3 C_3.
 \end{equation}

We now choose three subsets of the rows, columns and heights
 of $\cA$
 \begin{equation}\label{defIJK&card}
 I\subset \an{\ell_1},\ \#I=p,\quad J\subset \an{\ell_2},\ \#J=q,\quad
 K\subset\an{\ell_3},\ \#K=r.
 \end{equation}
 Let
 \begin{eqnarray}
 &&\cC_1=\cA_{\an{\ell_1},J,K}:=[a_{i,j,k}]\in\R^{\ell_1\times q\times r}, \; i\in
 \an{\ell_1},j\in J, k\in K, \nonumber\\
 &&\cC_2=\cA_{I,\an{\ell_2},K}:=[a_{i,j,k}]\in\R^{p\times \ell_2\times r},\;
 i\in I, j\in \an{\ell_2}, k\in K,
 \label{cchoice}\\
 &&\cC_3=\cA_{I,J,\an{\ell_3}}:=[a_{i,j,k}]\in\R^{p\times q\times \ell_3},\;
 i\in I, j\in J, k\in \an{\ell_3},\nonumber\\
 &&\cS=(\an{\ell_1}\times J\times K) \cup (I\times \an{\ell_2}\times K)\cup
 (I\times J\times \an{\ell_3}).\nonumber
 \end{eqnarray}
 We define $\cU_b$ and $\cU_{\opt }$ as in \cite{FMMN08}.
 \begin{eqnarray}\label{minprobtenb}
 \cU_b=\arg\min_{\cU\in\R^{m_1\times m_2\times m_3}}
 \sum_{(i,j,k)\in\an{l_1}\times\an{l_2}\times\an{l_3}} (a_{i,j,k}-\left(
 \cU\times_1 C_1 \times _2 C_2 \times _3 C_3
 \right)_{i,j,k})^2,\\
 \label{minprobtenopt}
 \cU_{\opt}=\arg\min_{\cU\in\R^{m_1\times m_2\times m_3}}
 \sum_{(i,j,k)\in\cS} (a_{i,j,k}-\left(
 \cU\times_1 C_1 \times _2 C_2 \times _3 C_3
 \right)_{i,j,k})^2.
 \end{eqnarray}

 Instead of computing $\cU_{\opt}$ we do the following
 approximations, as suggested in \cite{FMMN08} for the case $q=p,
 r=p^2$.  Unfold the tensor $\cA=[a_{i,j,k}]$ in the direction $3$ to obtain the matrix
 $E=[e_{s,k}]\in\R^{(\ell_1\cdot \ell_2)\times \ell_3}$.  So $e_{s,k}=a_{i,j,k}$
 for the corresponding pair of indices $(i,j)\in
 \an{\ell_1}\times\an{\ell_2}$.  Then the set of indices $(i,j)\in
 I\times J$ corresponds to the set of indices
 $L\subset\an{\ell_1\cdot\ell_2}$,  where $\#L=pq$.
 Denote by $ E_{L,K}$  the submatrix of $E$ which has row indices in $L$
 and column indices in $K$.  Let $E_{L,K}^{\dagger}\in
 \R^{r\times (pq)}$ be the Moore-Penrose inverse of $E_{L,K}$.
 As in \cite{FMMN08} we approximate the tensor $\cA$ by
 \begin{equation}\label{3tensdec}
 \cA_{\an{\ell_1},\an{\ell_2},K}E_{L,K}^{\dagger} \cA_{I,J,\an{\ell_3}}.
 \end{equation}
 For each $k\in K$ consider the matrix
 $$F_k:=\cA_{\an{\ell_1},\an{\ell_2},k}=[a_{i,j,k}]_{i,j=1}^{\ell_1,\ell_2}
\in \R^{\ell_1\times \ell_2}.$$
 Next we approximate $F_k$ by
 $G_k:=(F_k)_{\an{\ell_1},J}(F_k)_{I,J}^{\dagger} (F_k)_{I,\an{\ell_2}}$.
 As in \cite{FMMN08} we try several random choices of $I,J,K$
 with the cardinalities $p,q,r$ respectively, with the best
 preset conditions numbers for the matrices $E_{L,K}$ and
 $(F_k)_{I,J}$ for $k\in K$.

 Equivalently, we have that
 \begin{equation}\label{3tensdec1}
 \cA_{\an{\ell_1},J,k}\cA_{I,J,k}^{\dagger} \cA_{I,\an{\ell_2},k},
 \end{equation}
 is an approximation of $\cA_{\an{\ell_1},\an{\ell_2},k}$.
 Replacing $\cA_{\an{\ell_1},\an{\ell_2},k}$ appearing in (\ref{3tensdec})
 with the expression that appears in (\ref{3tensdec1}),
 we obtain the approximation $\cB$ of the form
 (\ref{tucrep2}).

 \bibliographystyle{plain}

\begin{thebibliography}{MMM}

 \bibitem{Fr08} S. Friedland, On the generic rank of 3-tensors,
 arXiv: 0805.3777v1.

 \bibitem{FriKNZ06} S. Friedland, M. Kaveh, A. Niknejad and H. Zare,
 Fast Monte-Carlo low rank approximations for matrices,
 \emph{Proc. IEEE Conference SoSE}, Los Angeles, 2006, 218-223.

 \bibitem{FMMN08} S. Friedland, V. Mehrmann, A. Miedlar and
 M. Nkengla, Fast low rank approximations of matrices and
 tensors, \emph{submitted}, www.matheon.de/preprints/4903.

 \bibitem{GolV96} G.H. Golub and  C.F. Van Loan, {\it Matrix
 Computation}, John Hopkins Univ. Press, 3rd Ed., 1996.

 \bibitem{GorT01} S.A. Goreinov and E.E. Tyrtyshnikov, The
 maximum-volume concept in approximation by low-rank matrices,
 \emph{Contemporary Mathematics} 280 (2001), 47-51.

 \bibitem{Har92} J. Harris, \emph{Algebraic Geometry},
 A First Course, 1992, Springer, New York.

 \bibitem{Kho} B.N. Khoromskij,
Methods of Tensor Approximation for Multidimensional
 Operators and Functions with Applications, \emph{Lecture at the workschop},
 Schnelle L\"oser f\"ur partielle Differentialgleichungen,
 Oberwolfach, 18.-23.05, 2008.

 \bibitem{deLdV00} L. de Lathauwer, B. de Moor and J. Vandewalle, A
 Multilinear singular value decomposition, \emph{SIAM J. Matrix
 Anal. Appl.} 21 (2000), 1253-1278.

 \bibitem{Lim07}
 V. de Silva and L.-H. Lim, "Tensor rank and the ill-posedness
 of the best low-rank approximation problem," to appear in SIAM
 Journal on Matrix Analysis and Applications.

 \bibitem{Tuc66} L.R. Tucker. Some mathematical notes on
 three-mode factor analysis, \emph{Psychometrika} 31 (1966),
 279 -- 311.









 \end{thebibliography}

\end{document}